\documentclass{amsart}

\theoremstyle{plain}
\newtheorem{thm}{Theorem}
\newtheorem{lem}[thm]{Lemma}

\theoremstyle{definition}
\newtheorem{defn}[thm]{Definition}
\newtheorem{ex}[thm]{Example}

\numberwithin{thm}{section}

\newcommand{\Lie}{\operatorname{Lie}}

\newcommand{\g}{\mathfrak{g}}
\newcommand{\frakb}{\mathfrak{b}}

\newcommand{\slt}{\mathfrak{sl}(3)}

\newcommand{\upg}{U_p(\g)}
\newcommand{\upigi}{U_{p_i}(\g_i)}
\newcommand{\uplgl}{U_{p_1}(\g_1)}
\newcommand{\uptgt}{U_{p_2}(\g_2)}
\newcommand{\dpigi}{D_{p_i}(\g_i)}
\newcommand{\dpg}{D_p(\g)}
\newcommand{\upgo}{U_p(\g,L)}

\newcommand{\upbp}{U_p(\frakb^+)}
\newcommand{\upbm}{U_p(\frakb^-)}
\newcommand{\upbpm}{U_p(\frakb^\pm)}
\newcommand{\upbmp}{U_p(\frakb^\mp)}
\newcommand{\upbip}{U_p(\frakb^+_i)}
\newcommand{\upblp}{U_p(\frakb^+_1)}
\newcommand{\upblm}{U_p(\frakb^-_1)}
\newcommand{\upbtp}{U_p(\frakb^+_2)}
\newcommand{\upbim}{U_p(\frakb^-_i)}
\newcommand{\uqblm}{U_q(\frakb^-_1)}
\newcommand{\uqblp}{U_q(\frakb^+_1)}
\newcommand{\upbipm}{U_p(\frakb^\pm_i)}
\newcommand{\upbimp}{U_p(\frakb^\mp_i)}

\newcommand{\upbimpl}{U_p(\frakb^\mp_i,L_i)}
\newcommand{\upbiml}{U_p(\frakb^-_i,L_i)}
\newcommand{\upblpl}{U_p(\frakb^+_1,L_1)}
\newcommand{\upblml}{U_p(\frakb^-_1,L_1)}
\newcommand{\upblpm}{U_p(\frakb^\pm_1)}
\newcommand{\upbtpm}{U_p(\frakb^\pm_2)}

\newcommand{\upibipm}{U_{p_i}(\frakb^\pm_i)}
\newcommand{\upbpl}{U_p(\frakb^+,L)}
\newcommand{\upbml}{U_p(\frakb^-,L)}

\newcommand{\catg}{\mathcal{C}(\Gamma)}

\newcommand{\C}{{\mathbb{C}}}
\newcommand{\Q}{{\mathbb{Q}}}
\newcommand{\Z}{{\mathbb{Z}}}
\newcommand{\vpi}{\varpi}

\newcommand{\cqg}{\C_q[G]}
\newcommand{\cpg}{\C_p[G]}
\newcommand{\cptg}{\C_{p,\tau}[G]}
\newcommand{\cqtg}{\C_{q,\tau}[G]}
\newcommand{\cpbipm}{\C_p[B_i^\pm]}
\newcommand{\cpblm}{\C_p[B_1^-]}
\newcommand{\cpbtp}{\C_p[B_2^+]}
\newcommand{\cpbip}{\C_p[B_i^+]}

\begin{document}


\title[Nonstandard quantum groups]{Nonstandard quantum groups
associated to certain Belavin-Drinfeld triples}
\author{Timothy J. Hodges}
\address{University of Cincinnati, Cincinnati, OH 45221-0025,
U.S.A.}
\email{timothy.hodges@uc.edu}
\thanks{The author was supported in part by NSF Grant DMS 9501484}

\begin{abstract}

     A construction is given of a family of non-standard
quantizations of the algebra of functions on a connected complex
semi-simple algebraic group. For each ``disjoint'' triple in the
sense of Belavin and Drinfeld, a 2-cocycle is constructed on
certain multi-parameter quantum groups. The new non-standard
quantum groups are the Hopf algebras obtained by twisting the known quantum groups
by these 2-cocycles. In particular,  the Cremmer-Gervais quantization of $SL(3)$ can be constructed in this way.

\end{abstract}

\maketitle

\section{Introduction}

	In \cite{EK1}, Etingof and Kazhdan proved  that any Lie bialgebra can be quantized. From this they deduced in \cite{EK2} that any Poisson algebraic group can be quantized.
Unfortunately the construction of Etingof and Kazhdan relies on a subtle existence theorem of Drinfeld and does not produce explicit descriptions of the algebras. It is therefore interesting to look for more explicit constructions of such quantum groups, in particular for important examples such as the Poisson groups given by Lie bialgebra structures that arise on semi-simple Lie algebras via skew-symmetric solutions of the classical Yang-Baxter equation \cite{BD}. Here we show how to construct some such non-standard quantum groups by twisting the standard ones by Hopf 2-cocycles. We give the construction for quantized function algebras in the case where the parameter $q$ is a complex number which is not a root of unity.  

	Let $\g$ be a complex Lie algebra. Then any solution of the classical Yang-Baxter equation defines a Lie bialgebra structure on $\g$. Of particular importance are the skew-symmetric solutions and these have been classified for a simple Lie algebra by Belavin and Drinfeld \cite{BD}. These are given by two parameters, one discrete (a `triple') and one continuous. Now let $G$ be a connected complex algebraic group with Lie algebra $\g$ and denote by $\cqg$ the standard quantized function algebra where $q \in \C$ is not a root of unity. Each of the Lie bialgebra structures on $\g$ defines a Poisson structure on $G$ and to each of these there should be a corresponding nonstandard version of $\cqg$.  How to deal with the continuous parameters is reasonably well understood (see for example \cite{HLT}). They correspond to twists of the quantum group by particularly elementary kind of 2-cocycle  or gauge transformation. We conjecture that for the discrete parameter there are more subtle kinds of 2-cocycles that can be constructed from these triples. We prove this conjecture here in a special case - when the triple is `disjoint'. 

	There are, of course, other ways of constructing non-standard quantum groups. In particular, given a solution $R$ of the quantum Yang-Baxter equation, one can use the technique of Fadeev, Reshetikhin and Takhtadjan to create a bialgebra $A(R)$ from which in certain situations one can pass to a Hopf algebra by taking a suitable quotient. This technique is used by Cremmer and Gervais \cite{CG} and Fronsdal and Galindo \cite{FG} to construct some interesting families of non-standard quantum groups. However the cocycle approach is far more powerful. There are none of the complications of passing from $A(R)$ to the Hopf algebra which become increasingly difficult as one gets away from the case $G=SL(n)$. The techniques below can, for instance, be applied to construct a number of non-standard versions of $\C_q[E_6]$. Moreover the representation (comodule) theory of a twisted quantum group is the same as that of the original quantum group and much of the theory of quantum $G$-spaces carries over similarly using routine arguments.

	Our construction covers many, but not all, of the known examples of explicit quantization of such Poisson groups. The examples constructed by Fronsdal and Galindo all correspond to disjoint triples and can therefore be approached from our point of view. The Hopf algebras underlying the quantum Lorentz group \cite{PW} and their generalizations the complex quantum groups \cite{CEJS} are quantizations corresponding to the natural triple on $\g \times \g$. This was pointed out by Majid in \cite{Ma}, though the point of view here is a little different. Of the Cremmer-Gervais quantizations of $SL(n)$ only the case $SL(3)$ is covered by our techniques. However, even in this case our construction provides interesting new information. For example, it shows that  the category of comodules is equivalent to the category of comodules over the standard quantum group $\C_q[SL(3)]$. 

	The construction outlined here is particularly simple in the situation where the triple is completely disjoint in the sense described after Theorem \ref{main}. In fact, in this case the cocyle can be constructed using a different (but essentially equivalent) procedure \cite{H2}. In the case where the triple is not completely disjoint, the idea behind the construction is still fairly simple but we are obliged to work with multi-parameter quantum groups and this complicates somewhat the technical details.

	This result was presented by the author at the JSRC meeting on {\em Quantization} at Mt Holyoke College. At this meeting, Etingof pointed out to the author that he had observed an analogous result for quantized universal enveloping algebras.


\section{Algebraic multi-parameter QUE's}

     Let $\g$ be a complex semi-simple Lie algebra  and let
$\Phi$ be a root system of $\g$ with respect to a fixed Cartan
subalgebra. Let $W$ be the Weyl group of $\Phi$. Let
$(\quad,\quad)$ be a $W$-invariant scalar product on the vector
space $\Q\Phi$ generated by $\Phi$ over the rationals. Fix a base
$\Pi$ of $\Phi$ and denote  the fundamental weights by
$\vpi_\beta$ for $\beta \in \Pi$.  Denote by $\Lambda$ the
weight lattice. Let $G$ be a connected complex algebraic group
with $\Lie G = \g$ and let $\Gamma$ be the sublattice of $\Lambda$ consisting of weights of representations of $\g$ induced from representations of $G$.

     Let $u: \Q\Phi \otimes \Q \Phi \to \Q$ be an alternating
bilinear form. Define forms $u_\pm$ by
$$
     u_\pm(\lambda,\mu) = u(\lambda,\mu) \pm (\lambda,\mu).
$$
Notice that the $u_\pm$ are non-degenerate since
$u_\pm(\lambda,\lambda) =(\lambda,\lambda)$ for all $\lambda \in
\Q \Phi$. Notice also that
$$
u_\pm(\lambda,\mu) = -u_\mp(\mu,\lambda).
$$
Let $\phi:\Q\Phi \to \Q \Phi$ be the map given by
$u(\lambda,\mu)=(\phi \lambda,\mu)$ and set $\phi_\pm = \phi \pm
1$ so that $u_\pm(\lambda,\mu)=(\phi_\pm \lambda,\mu)$.
Notice that $\phi_\pm$ are isomorphisms of $\Q\Phi$. Denote by
${\;\;}{\tilde{}}{\;\;}$ the isomorphism given by
$$
\tilde{\lambda} =-\phi_+^{-1}\phi_-(\lambda).
$$

     Choose $\hbar \in \C$ and set $q = \exp (- \hbar/2)$ and
more generally, for any $x \in \C$ set
$$
q^x = \exp (-x \hbar/2).
$$
Let $q_\alpha = q^{(\alpha,\alpha)/2}$.
Define $p,p_\pm : \Q\Phi \to \C$ by
$$p(\lambda,\mu) = q^{\frac{1}{2}u(\lambda,\mu)}, \quad
     p_\pm(\lambda,\mu) = q^{u_\pm(\lambda,\mu)}
$$
so that
$$p_\pm(\lambda,\mu) = p^2(\lambda,\mu)q^{\pm(\lambda,\mu)}
$$
Then $p$ is an antisymmetric bicharacter on $\Q\Phi$. Note also
that 
$$p_+(\tilde{\lambda},\mu) =p_-(\lambda,\mu)^{-1}$$
for all $\lambda, \mu \in \Q\Phi$.

     Given such a $p$ we may define a multiparametric version of
the quantized universal enveloping algebra. (The following construction is essentially that given in \cite{HLT} except that there we constructed the Drinfeld double rather than the quantized universal enveloping algebra. The definition we give here is the Drinfeld double constructed in \cite{HLT} factored out by the radical of the pairing between it and the quantized algebra of functions $\cpg$.) As in the one
parameter case, there are different versions of $\upg$
depending on a choice of lattice $\Omega$ between $\Z\Phi$ and
$\Lambda$.
We shall assume for now that $\Omega=\Z \Phi$ but the following
constructions apply equally well for any such $\Omega$. Set
$\hat{\Omega} = \Omega +
\tilde{\Omega}$. Define $\upg$ to be the algebra generated by
elements $E_\alpha$ and $F_\alpha$ for $\alpha \in \Pi$ and
$K_\lambda$ for $\lambda \in \hat{\Omega}$ subject to the
relations
$$
     K_\lambda K_\mu = K_{\lambda + \mu}
$$
$$
     K_\lambda E_\beta K_{-\lambda} = p_+(\lambda,\beta) E_\beta
$$
$$
     K_\lambda F_\beta K_{-\lambda} = p_+(\lambda,\beta)^{-1}
F_\beta
$$   
$$
     E_\alpha F_\beta - F_\beta E_\alpha = \delta_{\alpha
\beta}\frac{K_{\tilde{\alpha}}-K_{\alpha}^{-1}}{q_\alpha-q_\alpha
^{-1}}
$$
and the multiparameter quantum Serre relations,
$$\sum_{k = 0}^{1-a_{\alpha \beta}} (-1)^{k}
\left[ \begin{smallmatrix} 1-a_{\alpha \beta} \\ k \end{smallmatrix} \right]_{\alpha}
p(\alpha,\beta)^{-2k} E_\alpha^{1-a_{\alpha \beta}-k}E_\beta E_\alpha^k = 0, \mbox{ if } i \neq j $$
$$\sum_{k = 0}^{1-a_{\alpha \beta}} (-1)^{k}
\left[ \begin{smallmatrix} 1-a_{\alpha \beta} \\ k \end{smallmatrix} \right]_{\alpha}
p(\alpha,\beta)^{2k} F_\alpha^{1-a_{\alpha \beta}-k}F_\beta F_\alpha^k = 0,
\mbox{ if } i \neq j$$
where $a_{\alpha \beta}= 2(\alpha,\beta)/(\alpha,\alpha)$ and $\left[ \begin{smallmatrix} n \\ k \end{smallmatrix} \right]_{\alpha}$ is the usual quantum binomial coefficient.
The algebra $\upg$  has a Hopf algebra structure. The comultiplication and
counit are given by
$$
\Delta(E_\alpha) = E_\alpha \otimes 1 + K_{\tilde{\alpha}}
\otimes E_\alpha, \quad \epsilon(E_\alpha) = 0
$$
$$\Delta(F_\alpha) = F_\alpha \otimes K_{-\alpha} + 1 \otimes
F_\alpha, \quad \epsilon(F_\alpha)=0
$$
for all $\alpha \in \Pi$ and
$$
\Delta(K_\lambda) = K_\lambda \otimes K_\lambda, \quad
\epsilon(K_\lambda)=1
$$
for all $\lambda \in \hat{\Omega}$. The antipode $S$ is given by 
$$
S(E_\alpha) = -K_{-\tilde{\alpha}}E_\alpha, \quad S(F_\alpha) =
-F_\alpha
K_\alpha, \quad S(K_\lambda) = K_{-\lambda}.
$$

     For all $\nu \in \Omega$ define
$$
     U_\nu = \{ u \in \upg \mid K_\lambda u K_\lambda^{-1} =
          p_+(\lambda, \nu) u \}
$$
This defines an $\Omega$-grading on $\upg$.
 Inside $\upg$ we have two graded Hopf subalgebras,
$$
\upbp = \C\langle E_\alpha, K_\lambda \mid \alpha \in \Pi,
\lambda \in \tilde{\Omega} \rangle
$$
and
$$
\upbm = \C\langle F_\alpha, K_\lambda \mid \alpha \in \Pi,
\lambda \in \Omega \rangle
$$
Inside $\upbp$ are the graded subalgebras,
$$
U^+ = \C\langle E_\alpha \mid \alpha \in \Pi \rangle
\quad \text{and} \quad
\tilde{U}^0 = \C \langle K_\lambda \mid \lambda \in \tilde{\Omega} \rangle
$$ and inside $\upbm$ are the algebras
$$
U^- = \C\langle F_\alpha \mid \alpha \in \Pi \rangle
\quad \text{and} \quad
U^0 = \C \langle K_\lambda \mid \lambda \in \Omega \rangle
$$
Clearly $\upbp = \tilde{U}^0 U^+$ and $\upbm = U^0 U^-$.

     There is a nondegenerate skew Hopf pairing 
$$
\langle \; \mid \; \rangle: \upbp \otimes \upbm \to \C
$$
such that for all $x \in U^+$, $y \in U^-$, $\lambda, \mu \in \Omega$,
$$
     \langle xK_{\tilde{\lambda}} \mid y K_\mu \rangle = p_-(\lambda, \mu)
\langle x \mid y \rangle.
$$
Denote by $\dpg$ the associated Drinfeld double $\upbp
\bowtie \upbm$. There is a natural epimorphism of Hopf algebras
from $\dpg \to \upg$ given by $x \otimes y \mapsto xy$.

     For any $\upg$-module $M$ and any $\lambda \in \Lambda$,
set
$$
M_\lambda = \{ m \in M \mid K_\mu m = p_+(\mu,\lambda) m \text{
for all } \mu \in \hat{\Omega} \}.
$$
The elements of $M_\lambda$ are called weight vectors of weight
$\lambda$ and $M_\lambda$ is called the $\lambda$-weight space of
$M$. Recall that $\Gamma$ is the sublattice of $\Lambda$ of weights of representations of $G$. Denote by $\catg$ the subcategory of finite
dimensional $\upg$-modules $M$ such that
$$
     M= \bigoplus_{\lambda \in \Gamma} M_\lambda
$$
This category is a braided monoidal category  and is equivalent,
as a braided
monoidal category to the analogous category of modules over the
one parameter quantized universal enveloping algebra.
The restricted dual of $\upg$ with respect to $\catg$ is
denoted by $\cpg$.  In fact $\cpg$ is a cocycle twist of the one
parameter quantum group $\cqg$ with respect to a cocycle defined
using $p$ \cite{HLT}. 

     We will briefly need a slightly more general form of
$\upg$. Let $L$ now be any subgroup between $\Lambda$ and
$\Q \Phi$. Then we may define an algebra $\upgo$ exactly as
above. This
becomes a Hopf algebra for the analogous definitions of $S$,
$\Delta$ and $\epsilon$. 
The skew
pairing $\langle \quad,\quad \rangle$ extends to a nondegenerate
skew pairing on the more general versions of $\upbpl$ and
$\upbml$ which again satisfies
$$
     \langle xK_{\tilde{\lambda}} \mid y K_\mu \rangle = p_-(\lambda, \mu)
\langle x \mid y \rangle.
$$
for all $x \in U^+$, $y \in U^-$, $\lambda, \mu \in L$.

\section{Constructing 2-cocycles from Belavin-Drinfeld triples}

     A triple as defined by Belavin and Drinfeld, on a root
system $\Phi$ with base $\Pi$ is a triple $(\tau, \Pi_1,\Pi_2)$
where $\Pi_1,\Pi_2 \subset \Pi$ and $\tau: \Pi_1 \to \Pi_2$ is a
bijection satisfying,
\begin{enumerate}
\item $(\tau\alpha, \tau \beta) = (\alpha, \beta)$ for all
$\alpha, \beta \in \Pi_1$.
\item for all $\alpha \in \Pi_1$, there exists a $k$ such that
$\tau^k\alpha \not \in \Pi_1$
\end{enumerate}
We shall only be concerned here with the special case when $\Pi_1
\cap \Pi_2 = \emptyset$. We shall refer to such triples as {\em
disjoint triples}. 

\begin{defn}
An alternating bilinear form $u:\Q\Phi \otimes \Q \Phi \to \Q$ is
said to be compatible with a disjoint triple $(\tau,
\Pi_1,\Pi_2)$ if
\begin{enumerate}
\item $u(\tau\alpha, \tau \beta) = u(\alpha, \beta)$ for all
$\alpha, \beta \in \Pi_1$.
\item $u_+(\alpha, \tau\beta) = 0$ for all $\alpha, \beta \in
\Pi_1$.
\end{enumerate}
\end{defn}

Henceforth we shall fix a disjoint triple $(\tau, \Pi_1,\Pi_2)$
and a compatible alternating bilinear form $u$. We then proceed
to construct from $\tau$, a 2-cocycle on the quantum group $\cpg$
associated to $u$. 
Recall that a 2-cocycle on a Hopf algebra $A$ is an
invertible pairing $\sigma : A \otimes A \to k$ such that for all
$x$, $y$ and $z$ in $A$,
$$\sum \sigma(x_{(1)},y_{(1)}) \sigma(x_{(2)}y_{(2)},z) = \sum
\sigma(y_{(1)},z_{(1)})
     \sigma(x,y_{(2)}z_{(2)}) $$
and $\sigma(1,1) = 1$.
Given a 2-cocycle $\sigma$ on a Hopf algebra, one can twist
the multiplication to get a new Hopf algebra $A_\sigma$. The new
multiplication is given by
$$x \cdot y = \sum \sigma(x_{(1)},y_{(1)}) x_{(2)} y_{(2)}
\sigma^{-1}(x_{(3)},y_{(3)}). $$
See \cite{Mb} or \cite{DT} for further details.

     Associated to each $\Pi_i$ we have a root subsystem $\Phi_i
= \Phi \cap \Z\Pi_i$, a semi-simple Lie subalgebra $\g_i \subset
\g$ and a Weyl group $W_i$. The form $(\quad,\quad)$ restricts to
a $W_i$-invariant form on $\Q\Phi_i$. Similarly, $u$ restricts to
an alternating form $u_i$ on $\Q\Phi_i$. We thus have an
associated multiparameter quantized universal enveloping algebra
$\upigi$. The definition of a triple and the first compatibility
condition imply that $\tau$ induces an isomorphism of Hopf
algebras between $\uplgl$ and $\uptgt$. In general there does not
exist a homorphism from $\upigi$ to $\upg$. However the Hopf
subalgebras  $\upibipm$ do embed into $\upg$. Their  images are
the subalgebras $\upbipm$ where
$$
\upbip = \C \langle E_\alpha, K_{\tilde{\alpha}}^{\pm 1} \mid \alpha
\in \Pi_i
\rangle,
$$
$$
\upbim = \C \langle F_\alpha, K_\alpha^{\pm 1} \mid \alpha \in \Pi_i
\rangle.
$$
These maps can be combined to give a map from $\dpigi$ to $\upg$
which may or may not factor through to $\upigi$. Combining these
observations yields the following.

\begin{lem}
     There are Hopf algebra isomorphisms $\psi^\pm : \upblpm \to
\upbtpm$ given by $\psi^+(E_\alpha) = E_{\tau(\alpha)}$, 
$\psi^+(K_{\tilde{\lambda}}) =
K_{\widetilde{\tau(\lambda)}}$ and $\psi^-
(F_\alpha) = F_{\tau(\alpha)}$, $\psi^-(K_\lambda) =
K_{\tau(\lambda)}$.
\end{lem}

     The key to our construction is the following observation.

\begin{lem} The map 
$$
\phi: \upblm \otimes \upbtp \to \upg
$$
given by $\phi(u \otimes u') =uu'$ is a homomorphism of Hopf
algebras.
\end{lem}

\begin{proof}
     It suffices to notice that the generators of $\upblm$ and
$\upbtp$ commute. Let $\alpha \in \Pi_1$ and $\beta \in \Pi_2$.
Because $\Pi_1$ and $\Pi_2$ are disjoint,  $F_\alpha$ commutes
with $E_\beta$. On the other hand,
$$ K_\alpha E_\beta K_\alpha^{-1} = p_+(\alpha, \beta) E_\beta =
E_\beta$$
by the second part of the compatibility condition.  Similarly,
$$ K_{\tilde{\beta}} F_\alpha K_{\tilde{\beta}}^{-1} =
p_+(\tilde{\beta}, \alpha)^{-1} F_\alpha  = F_\alpha $$
since $ p_+(\tilde{\beta}, \alpha)= p_+( \alpha,\beta)=1$.
\end{proof}

	The map $\phi$ induces a map $\phi^*: \cpg \to (\upblm \otimes \upbtp)^\circ$. The next step is to identify carefully the image of $\phi^*$. The embeddings $\upbipm \to \upg$ yield maps $\rho_i^\pm :\cpg \to
\upbipm^\circ$. Denote the image of $\cpg$ in $\upbipm^\circ$ by
$\cpbipm$. Note that $\phi^* = (\rho_1^- \otimes \rho_2^+)\Delta$.

	Since $u_-$ is non-degenerate on $\Q\Phi_i$ we have projections $\pi_i^\pm : \Q \Phi \to \Q \Phi_i$ given by
$$
	u_-(\pi_i^+(\lambda),  \mu) = u_-( \lambda , \mu)  \quad\forall\; \mu \in \Q\Phi_i
$$
$$
	u_-(\mu, \pi_i^-(\lambda)) =  u_-(\mu , \lambda) \quad \forall\; \mu \in \Q\Phi_i.
$$
Set $L_i = \pi_i^-(\Z\Phi)$ and note that $\pi_i^+(\Z\Phi)=\tilde{L_i}$, where the tilde now refers to the automorphism of $\Q\Phi_i$ induced from the restriction of $u$ to $\Q\Phi_i$. Recall the definitions of $U_{p_i}(\frakb^\mp,L)$ from the end of the first section and denote these algebras by $\upbimpl$. Since there is a nondegenerate pairing between $\upbimpl$ and $\upbipm$, we can identify the former (after suitable reversal of the algebra or coalgebra structure) with a Hopf subalgebra of the dual of the latter.

\begin{lem}
      The algebras $\cpbipm$ coincide with the images of $\upbimpl$ respectively in $(\upbipm)^\circ$.
\end{lem}

\begin{proof}
The embeddings $\upbpm \to \upg$ yield maps $\rho^\pm :\cpg \to
\upbpm^\circ$. The Hopf pairing between $\upbp$ and $\upbm$ induces maps $\theta^\pm :\upbmp \to (\upbpm)^\circ$; $\theta^+$ is an antialgebra and coalgebra morphism, $\theta^-$ is an algebra and anti-coalgebra morphism. The images of $\rho^\pm$ coincide with those of $\theta^\pm$ \cite[Proposition 4.6]{HLT}. Similarly the pairing between $\upbip$ and $\upbim$ induces maps $\theta^\pm_i : \upbimp \to (\upbipm)^\circ$. Moreover for any $L$ between $\Omega$ and $\Q\Phi$, the maps $\theta^\pm_i$ extend to maps $\theta^\pm_i : \upbimpl \to (\upbipm)^\circ$.
       The maps $\rho_i^\pm$ factor through the maps $\rho^\pm : \cpg \to (\upbpm)^\circ$. Thus $\cpbipm$ is  the image of $\upbmp$ in $(\upbipm)^\circ$.  Recall that the pairing is given by 
$$
     \langle xK_{\tilde{\lambda}} \mid y K_\mu \rangle = p_-(\lambda, \mu)
\langle x \mid y \rangle.
$$
for all $x \in U^+$, $y \in U^-$, $\lambda, \mu \in \Omega$ and that for $\lambda, \mu \in \Z^+\Pi$, 
$$\langle U^+_\lambda, U^-_{-\mu} \rangle \neq 0 \text{ implies } \lambda = \mu.$$
The image of $U^\mp$ is thus easily seen to be $U^\mp \cap \upbimp$. For any $\lambda \in \Q\Phi_i$, $x \in U^+ \cap \upbip$ and $\mu \in \Q\Phi$,
$$
\langle x K_{\tilde{\lambda}} \mid K_\mu \rangle = p_-(\lambda, \mu) \epsilon(x)
= p_-(\lambda, \pi_i^-(\mu)) \epsilon(x)
= \langle x K_{\tilde{\lambda}} \mid K_{\pi_i^-(\mu)} \rangle.
$$
So the image of $K_\mu$ in $(\upbip)^\circ$ is $\theta^+_i(K_{\pi_i^-(\mu)})$. Thus $\cpbip = \theta^+_i(\upbiml)$.
\end{proof}

	As noted above,  there is a skew pairing  $\upblpl \otimes \upblml \to \C$. This induces, via $\theta^-_1\otimes (\theta^+_2\psi^-)$, a skew pairing on $\cpblm \otimes \cpbtp$. This skew pairing then induces a 2-cocycle on  $\cpblm \otimes \cpbtp$ given by
$$ \gamma_0(a \otimes b, c \otimes d) =
	\epsilon(a) \sigma_0(c,b) \epsilon(d)
$$
(see \cite{Mb} or \cite{DT}).
Finally let $\gamma: \cpg \otimes \cpg \to \C$ be the form induced on $\cpg$ via $\phi^*: \cpg \to \cpblm \otimes \cpbtp$.
Then $\gamma$ is a 2-cocycle on $\cpg$.  To summarize:

\begin{thm} \label{main}
     Let  $(\tau, \Pi_1,\Pi_2)$ be a disjoint triple, let $u$
be a compatible alternating form on $\Q \Phi$ and let $p$ be the
associated bicharacter. Then the map $\gamma$ is a 2-cocycle on
$\cpg$.
\end{thm}

	Let us denote the twisted Hopf algebra $(\cpg)_\gamma$ by $\cptg$. It is natural to ask whether $\cptg$ is different from $\cpg$. If $\tau$ is non-trivial, then $\cptg$ will never be isomorphic to $\cpg$ but we will not prove this here. Consider, instead, the special case where the triple $\tau$ is {\em completely disjoint} in the sense that $(\alpha, \beta)=0$ for all $\alpha \in \Pi_1$ and $\beta \in \Pi_2$. In this case the trivial pairing $u=0$ is compatible with $\tau$ and  $\phi^*$ maps surjectively from  $\cqg$ onto $\uqblp \otimes \uqblm$. Hence $\phi^*$ induces a surjective map from $\cqtg$ onto the double $\uqblp \bowtie \uqblm$. Clearly no such map exists for $\cqg$ since the irreducible finite dimensional representations of $\cqg$ are all one dimensional, whereas $\uqblp \bowtie \uqblm$ has irreducible representations of arbitrarily large dimension. Hence the algebra structure of $\cqtg$ is significantly different from that of $\cqg$.

\begin{ex}[Cremmer-Gervais $SL(3)$]
	Let $\g = \slt$ and $\Pi = \{\alpha_1,\alpha_2\}$. Let $\Pi_i = \{\alpha_i\}$ for $i = 1,2$ and let $\tau(\alpha_1) = \alpha_2$. The only $u$ compatible with $\tau$ is $u(\alpha_1, \alpha_2) = -1$. In this case $L_1 = \frac{1}{3} Z \alpha_1$, so that $\upblpl$ is the usual $\upblp$ extended by a `cube root' $K_1^{1/3}$. By looking in detail at the braiding on the category of $\C_{p,\tau}[SL(3)]$ modules one can show that  $\C_{p,\tau}[SL(3)]$ is precisely the quantization of $SL(3)$ associated with the R-matrix of Cremmer and Gervais \cite{CG,H1}. The map $\phi^*: \C_{p,\tau}[SL(3)] \to \uqblp \bowtie \uqblm$ coincides with an analogous map constructed in \cite[Theorem 5.5]{H1} for the Cremmer-Gervais quantizations of $SL(n)$.
\end{ex}

\begin{ex}
	More generally let $\g = \mathfrak{sl}(2n+1)$ and $\Pi = \{\alpha_1,\dots \alpha_{2n}\}$. Let $\Pi_1 = \{\alpha_1, \dots \alpha_n\}$ , $\Pi_2 = \{\alpha_{n+1}, \dots \alpha_{2n}\}$   and let $\tau(\alpha_i) = \alpha_{i+n}$. Then $\tau$ is again a disjoint triple and the above construction gives a nonstandard quantization of $SL(2n+1)$. In \cite{FG}, Fronsdal and Galindo constructed some $R$-matrices corresponding to these triples. Presumably the quantum groups associated to the $R$-matrices of Fronsdal and Galindo coincide with a quantum group of the form $\cptg$ for suitable choice of $p$.
\end{ex}

\begin{ex}[Double or complex quantum groups]
	If $\Pi$ is a base for $\g$, then $\Pi \sqcup \Pi$ is a base for $\g \times \g$. There is an obvious triple defined by letting $\Pi_1$ and $\Pi_2$ be the first and second copies of $\Pi$ respectively and letting $\tau$ be the natural bijection between them. Applying our construction yields the same cocycle as the one used in \cite{Ma} and \cite{H2} to construct the double or complex quantum groups of which the simplest is the quantum Lorentz group constructed by Podles and Woronowicz \cite{PW}.
\end{ex}

\end{document}